\documentclass[11pt]{article}
\usepackage{amsfonts, amssymb, graphicx, a4, epic, eepic, epsfig,
setspace, lscape}

\def\inst#1{$^{#1}$}

\newtheorem{theorem}{Theorem}[section]
\newtheorem{lemma}[theorem]{Lemma}
\newtheorem{proposition}[theorem]{Proposition}

\newtheorem{definition}[theorem]{Definition}
\newtheorem{corollary}[theorem]{Corollary}
\newtheorem{remark}[theorem]{Remark}

\newcommand{\cX}{{\cal X}}

\def \a {{\alpha}}
\def \b {{\beta}}
\def \e {{\varepsilon}}

\def \r {{\rho}}
\def \l {{\lambda}}
\def \m {{\mu}}

\def \s {{\sigma}}

\def \t {{\tau}}
\def \o {{\omega}}
\def \d {{\delta}}
\def \p {{\pi}}

\newcommand{\be}[1]{\begin{equation}\label{#1}}
\newcommand{\ee}{\end{equation}}

\newcommand{\bl}[1]{\begin{lemma}\label{#1}}
\newcommand{\el}{\end{lemma}}

\newcommand{\br}[1]{\begin{remark}\label{#1}}
\newcommand{\er}{\end{remark}}

\newcommand{\bt}[1]{\begin{theorem}\label{#1}}
\newcommand{\et}{\end{theorem}}

\newcommand{\bd}[1]{\begin{definition}\label{#1}}
\newcommand{\ed}{\end{definition}}

\newcommand{\bcl}[1]{\begin{claim}\label{#1}}
\newcommand{\ecl}{\end{claim}}

\newcommand{\bp}[1]{\begin{proposition}\label{#1}}
\newcommand{\ep}{\end{proposition}}

\newcommand{\bc}[1]{\begin{corollary}\label{#1}}
\newcommand{\ec}{\end{corollary}}

\newcommand{\bi}{\begin{itemize}}
\newcommand{\ei}{\end{itemize}}

\newcommand{\ben}{\begin{enumerate}}
\newcommand{\een}{\end{enumerate}}

\begin {document}

\title{Phase transitions for the 
cavity approach to the clique problem on random graphs
\thanks{Supported by GDRE 224 GREFI-MEFI and the European Research Council
through the ``Advanced Grant'' PTRELSS 228032}}

\author{
Alexandre Gaudilli\`ere\inst{1}\and
Benedetto Scoppola\inst{2}\and
Elisabetta Scoppola\inst{3}\and
Massimiliano Viale\inst{4}}

\date{}

\maketitle

\begin{center}
{\footnotesize
\vspace{0.3cm} \inst{1} LATP, Universit\'e de Provence, CNRS\\
39 rue F. Joliot-Curie - 13013 Marseille, France\\
\texttt{gaudilli@cmi.univ-mrs.fr}\\ 

\vspace{0.3cm} \inst{2} Dipartimento di Matematica, University of Rome
``Tor Vergata''\\
Via della Ricerca Scientifica - 00133 Rome, Italy\\
\texttt{scoppola@mat.uniroma2.it}\\

\vspace{0.3cm} \inst{3} Dipartimento di Matematica, University of Rome
``Roma Tre''\\
Largo San Murialdo, 1 - 00146 Rome, Italy\\
\texttt{scoppola@mat.uniroma3.it}\\ 

\vspace{0.3cm} \inst{4}  Dipartimento di Fisica, University of Rome
 Òla Sapienza"\\
P.le Aldo Moro, 2 00185 Rome, Italy\\
\texttt{massimiliano.viale@gmail.com}\\ }
\end{center}

\begin{abstract}

We give a rigorous proof of
 two phase transitions for 
 a disordered system designed to find large cliques inside Erd\"os random graphs.
 Such a system is 
associated with a conservative probabilistic cellular automaton
inspired by the cavity method originally introduced in spin glass theory.

\end{abstract}

\eject
\tableofcontents

\section{Introduction}

The largest clique problem (LCP), is the problem to find the largest complete subgraph of a given graph G. 
Let $G = (V, E)$ be a graph. A graph $g$ is a {\it subgraph of} $G$, $g \subset G$, if 
its vertex set $V(g) \subset V$ and its edges $E(g) \subset E$. A subgraph $g=(V(g),E(g))$ is {\it complete} if for any $i, j\in V(g)$ then $(i,j)\in E(g)$.
We will denote by ${\cal K}(G)$ the set of {\it complete subgraphs} or 
{\it cliques} of $G$ and by ${\cal M}axCl(G)$ the set of 
the largest cliques in $G$:
\be{maxcl}
{\cal M}axCl(G) := \{g \in {\cal K}(G) : |V(g)| = \max_{g'\in{\cal K}(G)} |V(g')| 
\ee
where $|B|$ denotes the cardinality of the set $B$.
We call {\it clique number} of the graph $G$, $\omega(G)$, 
the cardinality of the vertex set
of any largest clique in $G$, i.e., $\omega(G) = |V(g)|$ 
with $g \in {\cal M}axCl(G)$.
Solving the LCP for a given graph $G(V,E)$ implies
finding 
$\omega(G)$, and both problems are in fact in the same complexity class.
Note that we are not strictly following  the  definition in \cite{bollobas} since we are
using the term clique also for a non maximal complete
subgraph of $G$.
The LCP is one of the main example
of $NP$-hard problem. It has been proven (see e.g. [GJ] and references therein) 
to be polinomially equivalent to the $k$-satisfiability problem and it is
equivalent to many other well known difficult problems in combinatorial
optimization.

It is well known   that the LCP remains difficult 
also when restricted to typical instances of Erd\"os random graphs with finite fixed density $p$, i.e. of graphs with $n$ vertices, $|V|=n$,  in which each 
pair $(i,j)\in V\times V$ belongs to the edges set $E$ with independent 
probability $p$.  In particular it is well known that in such a random graph $G$
 it is very easy to find complete subgraphs 
$g\in{\cal K}(G)$ of size $|g|=-\frac{\log n}{\log p}$ but 
is difficult to find cliques that exceed this size, see below. 
The clique number
$\omega(G)$ is almost deterministic, in a sense that
will be stated more precisely below, and it is roughly speaking twice the size
of the cliques that are easy to find.
Recently some progress has been made in \cite{OE} to understand the intricated
landscape of the LCP for Erd\"os random graph, and therefore to show why at
the moment there are no algorithms able to find cliques of a size that exceeds
the easy one.

In a previous paper \cite{ISS} , in collaboration with Antonio Iovanella, two of us introduced an
algorithm to find cliques  inspired by the cavity method developed in the study of spin glasses.
This  Markov Chain Monte Carlo exhibits very good numerical performances, in the sense that,
although asymptotically it is not able to find cliques larger than the easy ones, for finite size
effects it find cliques very near to the largest also for quite large graphs.
The idea of the algorithm is the following: starting from a non feasible 
(i.e. non clique) configuration $\sigma$
of $k$ vertices of $V$, the algorithm chooses the next configuration assigning to each
new set $\s'$ of $k$ vertices of $V$ a probability proportional to $e^{-\beta [H_0(\sigma, \sigma')
+h(k-q(\sigma, \sigma'))]}$, 
where $\beta$
is a parameter called {\it inverse temperature}. The function $H_0$ is a non negative quantity defined
by the number of missing edges between the two configurations, i.e. the number of pairs
$(i,j)$ with $i\in\sigma$, $j\in \sigma'$ and $i\neq j$ such that $(i,j)\notin E(G)$.
The quantity $q(\sigma, \sigma')$ represents the overlap between $\sigma$ and $\sigma'$
and then $k-q(\sigma, \sigma')$ is the number of vertices in $\sigma'$ that are not
in $\sigma$. The transition probabilities depend therefore also on the positive parameter $h$.

The presence of $\beta H_0$ in the transition probabilities, when  $\beta$ is large,
makes very low the probability to reach new configurations
$\sigma'$ that are badly connected with $\sigma$,  while a large $h$ depresses
the configurations $\sigma'$ with many different vertices with respect to $\sigma$.

From a statistical mechanics point of view the dynamics above has various interesting features.
First, 
the dynamics is conservative, since it is defined
on the space of configurations with $k$ vertices,
moreover since the whole configuration can be renewed in a single
step the resulting MCMC can be considered a {\em canonical} (or {\em conservative})
probabilistic cellular automaton (PCA).
Rigorous results on canonical PCA are quite rare in the literature. 
Second, $H(\sigma, \sigma')=H_0+h(k-q)$ is in 
some sense the Hamiltonian of a disordered system of pair of configurations, and the 
combined action of $ H_0$ and $ h(k-q)$ 
makes the energy landscape quite complicate. 
Third, good numerical performances
stimulate a deeper understanding of the dynamics.

For this reasons
we decided to study in more detail the statistical mechanical system described by the chain 
in the case of random graphs. We prove several results. 
First of all it can be proved rigorously that for suitable values of $k$,
including the interesting case $k=\o(G)$,
 the annealed analysis corresponds to the quenched one.

Then it can be proved 
the existence in the plane $(\beta, h)$ of a nontrivial phase diagram.
More precisely, the system exhibits a first order phase transition while the pair $\beta, h$
crosses a line $h_c(\beta)$.
At $h>h_c(\beta)$ the phase is characterized by pairs of configurations $\s,\s'$ with
 $\sigma=\sigma'$ and with a given density of missing links in $\s$, depending on $\b$.
At $h<h_c(\beta)$ the phase is characterized by pairs of  disjoint configurations $\s,\s'$ with 
again a particular value for the density of missing links between $\s$ and $\s'$ depending on $\b$.

Moreover, in the region below the critical line $h_c(\beta)$, a second phase transition
is present, and again it has a transparent ``physical" interpretation: for temperature $T={1\over\beta}$
below a critical value $T_c$ the system tend to oscillate indefinitely between two fixed 
configurations $\sigma$ and $\sigma'$, while above $T$ the new configuration at each step is
typically different from the configurations previously visited by the system.

This detailed control on the features of the system is achieved by a careful evaluation
of the thermodynamics. 
In particular the proof of the existence of the phase transitions can be
performed in a relatively easy way, computing the annealed partition function of the system.
The self averaging of the system, i.e., the equivalence between quenched and annealed,
 is more complicate to prove rigorously, involving
the computation of the second moment of the partition function, and it is more a
brute-force computation.
We will present it in some detail in the paper, in an almost pedagogical way,
because, as far as we know, there are few cases in the literature where
a phase transition for 
 a disordered system can be controlled rigorously.
  Moreover the way we achieve this result, although based on 
classical argument like the saddle point method, has some technical details that are quite
interesting and may be useful also in different contexts.

Of course this analysis gives important information on the choice of the parameter
used in simulation, and in a following paper we will discuss its application to the study of
the convergence to equilibrium of the dynamics.

 To be more precise we need now some definition.


\subsection{Random graphs and the clique number}\label{RG}

In this section we fix definitions and notations on random graphs and we recall well known results
on the clique number.

For all  $p\in [0,1]$ consider the {\it probability space} 
given by an infinite sequence of independent Bernoulli variables of parameter $p$, i.e.,
$\o\in\Omega:=\{0,1\}^{\mathbb N},\; \o=(a_1,a_2,...,a_l,...)$ with $a_l\in\{0,1\}$, with $\s$-algebra
generated by $A_l^j:=\{\o:\; a_l=j\},\; j=0,1$ and with probability measure
$$
\mathbb P(\o:\;a_{i_1}=j_1,...a_{i_k}=j_k)=p_{j_1}...p_{j_k}\quad \hbox{ with } p_1=p,\; p_0=1-p.
$$
Given a set of vertices $V=\{1,...,n\}$ we associate to it the probability space $\Omega_n$ given by the first
${n\choose 2}$ Bernoulli variables in $\Omega$ describing  the  edges between vertices
in $V$, with
the obviuos ordering
$(1,2), (1,3), (2,3),....,(1,n), (2,n),..., (n-1,n)$. In this way we represent with $\Omega$ the
probability space usually denoted by ${\cal G}(\mathbb N,p)$, i.e., the infinite random graph.

For any $G\in{\cal G}(\mathbb N,p)$ and $n\in\mathbb N$  we denote by $G_n$ 
the subgraph of $G$ {\it spanned} by the set $V_n:=\{1,2,..., n\}$, i.e., the subgraph of $G$ containing
all the edges of $G$ that join two vertices in $V_n$. By definition $G_n$ is $\Omega_n$ measurable.
We will denote by $\mathbb P$ and $\mathbb E$  the probability and the mean value respectively, on this probability space.

The following well known result on the clique number can be found in \cite{bollobas}(Corollary 11.2, pg 286):

\bp{clique}
For a.e. $G\in {\cal G}(\mathbb N,p)$ there is a constant $m_0(G)$ such that if 
$n\ge m_0(G)$ then 
$$
\Big| \o(G_n)-2\log_b n+2\log_b\log_b n -2\log_b({e\over 2})-1\Big|<{3\over 2}
$$
with $b:={1\over p}$.
\ep
The main tool in  the proof of this Proposition is the study
 of the random variable $Y_r(n)$
 defined as the
number  of complete subgraphs of $G_n$ with $r$ vertices, i.e., the number of $r$-cliques
in $G_n$ with the second moment method . Indeed its mean value is given by:
$$
\mathbb E Y_r(n)={n\choose r}p^{{r\choose 2}}=:f(r,n)\simeq b^{r\log_b n-{r^2\over 2}}
$$
with $b:={1\over p}$.
The function $f(r,n)$, as a function of $r$, has its maximum in $r\simeq \log_b n$ and  drops rather suddenly below 1, by increasing
$r$, say at $r\sim 2\log_b n$. 
Moreover, again by an explicit calculation, $Y_r(n)$ satisfies the following inequality
$$
{varY_r(n)\over ( \mathbb E Y_r(n))^2}\le br^4 n^{-2}+2 (\mathbb E Y_r)^{-1}, 
$$
when $(1+\eta)\log_b n<r<3\log_b n$, for $\eta\in (0,1)$.
With the Borel-Cantelli lemma,  it is easy to show that, given
$\e\in (0,{1\over 2})$,  
for almost every graph $G\in {\cal G}(\mathbb N,p)$  there is a constant $m_0(G)$ such that if 
$n\ge m_0(G)$ and $n'_r\le n\le n_{r+1}$ then $\o(G_n)=r$, with
\be{nrn'r}
n_r:=\max\{n\in\mathbb N:\, f(r,n)\le r^{-(1+\e)}\}\quad
n'_r:=\min\{n\in\mathbb N:\, f(r,n)\ge r^{1+\e}\}.
\ee
Indeed  the size $k$ of the interesting cliques can be parametrized by a real $c\in (1,2]$
since the relation between $n$ and the size $k$ of the cliques that we want to study
is given by
$$
\ln n =k {\ln 1/p\over c}, \quad \hbox{with } c\in(1,2].
$$
As emerges in (\ref{nrn'r}) it's more efficient to use $k$ as parameter, instead of $n$, to study the asymptotic behavior for large graphs and so 
for any $\bar c>1$ we define 
\be{Sc}
{\cal S}_{\bar c}:=\{(n_k)_{k>0}:\; \lim_{k\to\infty}{\ln n_k\over k}={\ln 1/p\over \bar c}\}
\ee
This means that if we define $c_k=k{\ln 1/p\over\ln n_k}$ we consider sequences $n_k$
such that  $\lim_{k\to\infty}c_k=\bar c$. 
This is actually a particular asymptotic regime that could be generalized.

Let $Y$ be a {\it random function} on the probability space
$\Omega$  associating to a pair $({n,k})$ a random variable 
$Y(n,k)$ on $\Omega_n$ depending on $k$, for instance the number of $k$ cliques in $G_n$, considered before.
\bd{sa}
A  random function $Y$ on the probability space  $\Omega$  is called  {\it $\bar c$-asymptotically self averaging},
if the random variables
${Y(n_k,k)\over {\mathbb E}Y(n_k,k)}$ converge almost surely to $1$ for any $(n_k)_{k>0}$ in ${\cal S}_{\bar c}$
as $k\to\infty$, uniformly in ${\cal S}_{\bar c}$.  
\ed
This means that there exists $\tilde \Omega\subset \Omega$ with ${\mathbb P}(\tilde\Omega)=1$
such that for every $\o\in\tilde\Omega$ and any $(n_k)_{k>0}\in {\cal S}_{\bar c}$ the random variable
${Y(n_k,k)(\o)\over {\mathbb E}Y(n_k,k)}$ converges to $1$.

Note that, by the Borel-Cantelli lemma a sufficient condition for self-averaging is the following:
\be{selfa_var}
{var Y(n_k,k)\over  ({\mathbb E}Y(n_k,k))^2}<n_k^{-\a} e^{o(k)}
\ee
for some $\a>1$ with $o(k)$ uniform in $(n_k)\in {\cal S}_{\bar c}$.
Indeed for any $\e>0$ we have that 
$$P\Big(|{Y(n,k)\over {\mathbb E}Y(n,k)}-1|>\e\,\hbox{ for some }
n=n_k, \,(n_k)_{k>0}\in{\cal S}_{\bar c}\Big)\le 
$$
$$\le e^{{k\ln 1/p\over c}+o(k)}{1\over \e^2} var({Y(n_k,k)\over {\mathbb E}Y(n_k,k)})
\le {1\over \e^2} e^{k{\ln 1/ p\over c}(\a-1)+o(k)}
$$
is summable on $k$, since $(n_k)_{k>0}\in{\cal S}_{\bar c}$ implies that $c_k:=k{\ln 1/ p\over \ln n_k}$
converges to $\bar c$,  so that with at most finitely many exception on $k$, we have that
$|{Y(n_k,k)\over {\mathbb E}Y(n_k,k)}-1|<\e$ for any $(n_k)_{k>0}\in{\cal S}_{\bar c}$.

We also note that  if  $ Z$ is 
$\bar c$-asymptotically self averaging we have that $\ln Z(n_k,k)-\ln {\mathbb E}Z(n_k,k)$ converges
almost surely to $0$.  

\subsection{The cavity algorithm }\label{cliquealgo}

\noindent
Let $V=\{1,...,n\}$ and define for each unordered pair 
 in $V\times V$ 
\be{Jij}
J_{ij}=\cases{0&if $(i,j)\in E$\cr
1&if $(i,j)\notin E$\cr}
\ee

We consider the  space $\cX^{(n)}:= \{0,1\}^{\{1,...,n\}}$ of lattice gas
configurations on $V$ and we will denote by the same letter
a configuration $\s\in\cX^{(n)}$ and its support $\s\subseteq V$.
On this configuration space
$\cX$ we can consider an Ising Hamiltonian with an antiferromagnetic
interaction between non-neighbor sites:

\be{defham2}
H(\s):=   \sum_{i,j\in V,\; i\not=j}J_{ij}\s_i\s_j  -{h}\sum_{i\in V}\s_i
\ee
where  $h>0$.
It is immediate to prove that when  $h<2$  the minimal value of  $H(\s)$
is obtained on configurations with support on the vertices of a
maximum clique. 
In the case of a random graph $G$, i.e., when the interaction
variables $J_{ij}$ are i.i.d.r.v., the Hamiltonian (\ref{defham2})
is similar to the Hamiltonian of the Sherrington-Kirkpatrick(SK) model. The main differences
are that we use lattice gas variables instead of
spin variables and, more important,  the interaction 
is given by Bernoulli variables.

For each  $\s\in\cX^{(n)}$
 we define its {\it cavity field} (or {\it molecular field})
 as  the field  created in each site $i$ by
all the sites  in the configuration $\s$:
\be{hi} 
h_i(\s)= \sum_{j\not= i}J_{ij}\s_j+h(1 -\s_i) \qquad\forall
 i\in V.
\ee

We consider  the canonical case, i.e., for any integer $k<n$ we define
  the {\it canonical configuration space}
\be{Xk}
\cX_k^{(n)}:=\{\s\in\cX^{(n)}:\; \sum_{i\in V}\s_i=k\}
\ee

For each  pairs of configurations $\sigma,\sigma' \in \cX_k^{(n)} $ 
we can define the {\it pair hamiltonian}:
\be{Hss'}
H(\sigma,\sigma')=\sum_{i,j\in V, \, i\not=j}
J_{ij}\sigma_i\sigma'_j+h\sum_i(1-\s_i)\s'_i=\sum_ih_i(\s)\s'_i
\ee
This hamiltonian is non-negative and vanishes when
$\s=\s'$ and its support is a $k$-clique.

For every $\s,\s'\in\cX_k^{(n)}$ the {transition probabilities} of the {\it cavity algorithm} are given by:
\be{pss'}
P(\sigma,\sigma')=\frac{e^{-\beta  H(\sigma,\sigma')}}
{\sum_{\tau\in \cX_k} e^{-\beta H(\sigma,\tau)}}={e^{-\beta  H(\sigma,\sigma')}\over
{Z_\s}},
\ee
with
\be{Zs}
Z_\s=\sum_{\tau\in\cX_k} e^{-\beta H(\sigma,\tau)}.
\ee

By an immediate computation we can check that
 the detailed balance condition w.r.t. the {\it invariant measure} on $\cX_k^{(n)}$
\begin{equation}
\m(\sigma)=\frac
{\sum_{\tau\in\cX_k^{(n)}} e^{-\beta H(\sigma,\tau)}}
{\sum_{\tau,\sigma\in\cX_k^{(n)}} e^{-\beta H(\sigma,\tau)}}={Z_\s\over Z}
\label{Pi}
\end{equation}
is verified with  {\it partition function} $Z$:
\be{Z}
Z(n,k)=Z(\cX_k^{(n)})=\sum_{\s\in\cX_k^{(n)}}Z_\s=\sum_{\s,\t\in\cX_k^{(n)}}e^{-\b H(\s,\t)}.
\ee
We will denote by $\m(.)$ the mean w.r.t. this  stationary measure.
For large $\b$, this
stationary measure is exponentially concentrated on cliques.

Note that at each step all the sites are updated;  this dynamics could be considered 
a canonical  version of  probabilistic cellular automata (PCA).
Given a fixed configuration $\s$ 
 the probability measure on $\cX_k^{(n)}$  given by $\p_\s(.):=P(\s,.)$
can be considered in the frame of the {\it Fermi statistics}.
Indeed the cavity fields
$h_i(\s)$ have values $e_{l,r}=l+rh$ with
${l\in\{0,1,...,k\}}$ and $r\in\{0,1\}$. We define 
\be{Ij}
{\cal I}_{l,1}:=\{i\in
V:\, h_i(\s)=l+h\}\qquad l=0,...,k.
\ee 
By equation (\ref{Hss'}) we have
\be{Hen}
H(\s,\t)=\sum_ih_i(\sigma)\t_i=\sum_{l,r} e_{l,r} \sum_{i\in {\cal I}_{l,r}}\t_i=:\sum_{l,r}e_{l,r}n_{l,r}
\ee
where $n_{l,r}$ denotes the occupation of the level (or cell) $(l,r)$.
On the other hand each level  consists of $g_{l,r}:=|{\cal I}_{l,r}|$ subcells (or sublevels)
containing at most one particle since for every $i\in {\cal I}_{l,r}$ we have $\t_i\in\{0,1\}$.
This means that instead of  configurations  in $\cX_k^{(n)}$
we can consider the occupation numbers of the levels $\{n_{l,r}\}_{l=0,...,k,\,r=0,1}$.
This statistical system is called a {\it Fermi gas}, see \cite{AG} for more detail 
on sampling for the Fermi statistics, and thus
on the realization of this single step of the dynamics.

As far as the energy levels corresponding to sites not in $\s$, i.e., with $r=1$, are concerned,
 we have that their number of sublevels, $g_{l,1}=|\{i\in V:\; h_i(\s)=l+h\}|$, is ``almost deterministic",
 as discussed in \cite{ISS}. Indeed they follow a binomial law, and precise results can be found in
 Lemma \ref{nest} in Section \ref{S5}.

A final remark on probability measures can be useful.
The invariant measure $\m(\s)$ is not a Gibbs measure, as usual with PCA, 
but we can define a Gibbs measure on pairs of configurations, i.e., on
$\cX^{(n)}_k\times\cX^{(n)}_k$ as 
 $\m_2(\s,\s')={1\over Z} e^{-\b H(\s,\s')}$ with the same partition function $Z(n,k)$ given
 in (\ref{Z}). Actually the invariant measure $\m$ can be considered  the marginal of $\m_2$.
 The probability measure $\p_\s(.)={e^{-\b H(\s,.)}\over Z_\s}$ introduced above in the discussion on the Fermi statistics can be considered as the {\it conditioned measure}  on ${\cal X}^{(n)}_k$,
 since we have the relation:
\be{measures}
\m_2(\s,\t)=\m(\s)\p_\s(\t).
\ee


\subsection{Results}

The main results presented in this paper are summarized by the following:

\bt{T1}
For each $p\in(0,1)$ and  $\b\in(0,\infty]$

\bi
\item[i)] let $n$  and $k$ be integers  such that
$c:=k{\ln 1/p\over \ln n}>1$, 
defining $\tilde h:={h\over k}$,
there is a critical value of $\tilde h$ defined by
\be{tildehc}
\tilde h_c= {1\over\b}\Big({f(2\b)\over 2}-f(\b)+{\ln(1/p)\over c} \Big)
\ee
with $f(\b):=-\ln \big[p+(1-p)e^{-\b}\big]$ for which
\be{lnZ}
\ln ({\mathbb E}Z(\cX_k^{(n)}))= 
\cases{k^2 {\ln(1/p)\over c} +k-f(2\b){k(k-1)\over 2}-k\ln k+{ o}(k)
&if $\tilde h> \tilde h_c$\cr
2k^2 {\ln(1/p)\over c} +2k-\b \tilde hk^2-f(\b)k^2
-2k\ln k+{ o}(k)&if $\tilde h< \tilde h_c$\cr
}
\ee
\item[ii)]
if $c\in(1,2]$ 
\be{stvar}
{var Z(\cX_k^{(n)})\over \big({\mathbb E}Z(\cX_k^{(n)})\big)^2}\le n^{-2} e^{o(k)}
\ee
with $o(k)$ independent of $n$ and $c$,
so that the partition function $Z$ is $\bar c$-asymptotically self averaging for $\bar c\in (1,2]$;
\item[iii)]
 the line
$\tilde h_c$ corresponds to a first order phase transition, in particular the
phase with $\tilde h>\tilde h_c$ is characterized by configurations $\s$ with $H(\s,\s)\sim k^2 f'(2\b)$ and the phase with $\tilde h<\tilde h_c$ is characterized by pairs of disjoint
configurations $(\s,\s')$ with $H(\s,\s')\sim k^2f'(\b)+k^2\tilde h$;

\item[iv)] 
if $\bar c\in(1,2]$ and
for any $(n_k)_{k>0}\in{\cal S}_{\bar c}$ and  $\tilde h\not=\tilde h_{\bar c}$ define  the entropy  $S(n_k,k):= -\sum_{\s\in\cX_k^{(n_k)}}\m(\s)\ln(\m(\s))$, 
then ${1\over k^2}S(n_k,k)$ converges almost surely to the following non random
function:
\be{Sfin'}
\bar s=\cases{
 {\ln 1/p\over {\bar c}}-{f(2\b)\over 2}+\b f'(2\b) & in the parameter region  $(A):\;\tilde h>\tilde h_{\bar c}$\cr
2 {\ln 1/p\over {\bar c}}-f(\b)+\b f'(\b) & in the  region  $(B):\;\tilde h<\tilde h_{\bar c}$ and $\b>\b_{\bar c}$\cr
 {\ln 1/p\over {\bar c}}& in the  region  $(C):\;\tilde h<\tilde h_{\bar c}$ and $\b<\b_{\bar c}$\cr
}
\ee
where $\b_{\bar c}$ is a zero of the function
\be{Cb}
C(\b)={\ln 1/p\over {\bar c}}-f(\b)+\b f'(\b)
\ee
so that also $S(n_k,k)$ is ${\bar c}$-asymptotically self-averaging.
 The function $\bar s$ is discontinuous along the line
$\tilde h=\tilde h_{\bar c}$ and  has a  discontinuity in its first derivative at $T_{\bar c}={1\over \b_{\bar c}}$
corresponding to a ``low temperature phase transition".
The asymptotic value of the entropy in the phase $\tilde h<\tilde h_{\bar c}$ and $\b<\b_{\bar c}$ is maximal
since $|\cX_k^{(n_k)}|\asymp e^{k^2 {\ln 1/p\over c_k}}$.
\ei
\et
\begin{figure}
\begin{center}
\includegraphics[width=2in]{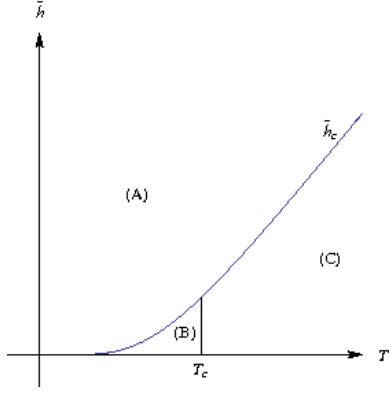}

\caption[]{Phase dygram for the cavity algorithm}
\end{center}      
\end{figure}

The phase diagram is summarized in Figure 1.

\br {r1}
We can write
\be{ZE}
 Z=\sum_{E} e^{-\b E} N_E
\ee
with $E$ running on all the possible values of the energy $H(\s,\t)$ and
$N_E$ being the number of pairs of configurations $\s,\t$ with  $H(\s,\t)=E$.
For $\b=\infty$ this implies that $Z=N_0$. Hence the self-averaging of $Z$
implies the self-averaging of the number of cliques of any size $k$ corresponding to $c\in(1,2]$.
This generalizes the Bollobas result quoted above.
\er
\br{r2}
Even though the relevant case for the clique problem is $c\in(1,2]$, 
for  $c>2$ we can prove (see Appendix~\ref{A2})
that, with $\bar\b_c<\infty$ the unique solution of 
$$
f(2\bar\b_c)-\frac{1}{2}f(4\bar\b_c)={1\over c}\ln 1/p,
$$
for all $\b<\bar\b_c$
we still have the estimate (\ref{stvar}).
Therefore quenched quantities behave like annealed ones for $\beta<\bar\beta_c$.
In addition we can prove the existence of a second value for the inverse temperature, say
$\hat\beta_{ c}  > \bar\beta_{ c}$ such that  for $\beta > \hat\beta_{ c}$ 
and $\beta > 2\hat\beta_{c}$, respectively for    $\tilde h >\tilde h_{ c}$ and  $\tilde h <\tilde h_{ c}$,
quenched quantities certainly differ from annealed ones.
Indeed if 
 $\hat\beta_{ c} > \bar\beta_{ c}$ is such that
$$
f(2\hat\b_{ c})-2\hat\beta_{ c}f'(4\hat\b_{ c})={1\over c}\ln 1/p,
$$ 
the estimated entropy for $\mu_2$ obtained from the annealed quantities
turns out to be asymptotically negative, 
i.e., for $(n_k)_{k>0}$ in ${\cal S}_{\bar c}$,
$$
\lim_{k\rightarrow +\infty}
\frac{1}{k^2}
\left(
  \ln{\mathbb E}[Z(n_k,k)] 
  - \beta\frac{\partial \ln{\mathbb E}[Z(n_k,k)]}{\partial \beta}
\right)
<0
$$
 for $\beta > \hat\beta_{\bar c}$ in the case $\tilde h >\tilde h_{\bar c}$
and $\beta > 2\hat\beta_{\bar c}$ in the case $\tilde h < \tilde h_{\bar c}$.
Then, for $\beta > \hat\beta_{\bar c}$ 
and $\beta > 2\hat\beta_{\bar c}$ respectively,
quenched quantities certainly differ from annealed ones.
We actually expect a third phase transition at these temperatures:
conversely, for $\beta < \hat\beta_{\bar c}$ 
and $\beta < 2\hat\beta_{\bar c}$ respectively,
quenched quantities should behave like quenched ones.


\er

{\bf Notation:}

\noindent
For notation convenience in what follows we adopt the following simplification:
given $\bar c>1$ and $(n_k)_{k>0}\in {\cal S}_{\bar c}$ we write $n=n_k$ and $c=c_k$


\section{The annealed partition function $\bar Z$}
\label{S2}
\bigskip

In the case of random graphs is not difficult to compute the annealed partition function: 
\be{Zann}
\bar Z:={\mathbb E}Z={\mathbb E}\Big[\sum_{\s,\t\in\cX_k}e^{-\b H(\s,\t)}\Big].
\ee

Let $I:=\s\cap\t$, and  $q$ be the overlap $q:=|I|$. We denote by $H_0(\s,\t)$
the first part of the pair hamiltonian, i.e., the pair hamiltonian evaluated for $h=0$:
\be{H0}
H_0(\s,\t)=\sum_{i,j\in V, \, i\not=j}
J_{ij}\sigma_i\sigma'_j= \sum_{\{i,j\}} J_{ij}(\s_i\t_j+\s_j\t_i)
\ee
The quantity $\s_i\t_j+\s_j\t_i$ takes values $0,1,2$ as given in Table 1,
where we denote by $S:=\s\backslash I$, $T:=\t\backslash I$ and $C:=(\s\cup \t)^c$.
 \begin{table}
\begin{center}
\begin{tabular}{|c|c|c|c|c|}
\hline
&S&I&T&C\\
\hline
S&0&1&1&0\\
\hline
I&1&2&1&0\\
\hline
T&1&1&0&0\\
\hline
C&0&0&0&0\\
\hline
\end{tabular}
\end{center}
\caption{Values of $\s_i\t_j+\s_j\t_i$}
\end{table}
So  $\s_i\t_j+\s_j\t_i=1$ for $\{i,j\}$ in the set of unordered  pairs
$ {\cal E}_1= (S\times I)\cup (T\times I)\cup (S\times T)$ and 
 $\s_i\t_j+\s_j\t_i=2$ for $\{i,j\}\in{\cal E}_2= I\times I$.
 By using the independence of the
random variables $J_{ij}$ we can conclude:
 $$
\bar Z =
\sum_{q=0}^k \sum_{\s,\t\in\cX_k:\; |I|=q}e^{-\b h (k-q)} \prod_{ \{i,j\}\in {\cal E}_1}
{\mathbb E} e^{-\b J_{ij}}
\prod_{ \{i,j\}\in {\cal E}_2}
{\mathbb E} e^{-2\b J_{ij}}.
$$
Since ${\mathbb E} e^{-\b J_{ij}}= e^{-f(\b)}$ and ${\mathbb E} e^{-2\b J_{ij}}=e^{-f(2\b)}$
with $f(\b)=-\ln\big[p+(1-p)e^{-\b})\big]$  and defining
\be{theta}
\Theta(q):= \ln \Big[ {n\choose 2k-q}{2k-q\choose q}{2(k-q)\choose k-q}\Big]
\ee
we have
\be{barz1}
\bar Z=\sum_{q=0}^k 
e^{\Theta(q)}
e^{-\b h(k-q)-f(\b)(2q(k-q)+(k-q)^2)-f(2\b){q(q-1)\over 2}}=\sum_{q=0}^k 
e^{\Theta(q)+\Phi(q)}
\ee
with
\be{phi}
\Phi(q):=-\b h(k-q)-f(\b)(k^2-q^2)-f(2\b){q(q-1)\over 2}
\ee
We collect in Appendix \ref{A0} the main properties of the function $f(\b)$.

As far as the entropic term $\Theta$ is concerned we can  use the Stirling formula $n!=({n\over e})^n n^{1\over 2}\sqrt{2\pi} e^{{\cal O}({1\over 12n})}$ 
to obtain the following asymptotic behaviors for 
$m$ and $k$ large with $m=o(n)$ and $g<k$:
\be{nk}
{n\choose m}=\exp \{m\ln n-m\ln m+m+{o}(m) \},\quad \ln (n-m)=\ln n-{m\over n}+{\cal O}({m^2\over n^2})
\ee
\be{kg}
{k\choose g}= \exp\{k\ln k-(k-g)\ln (k-g) -g\ln g+{ o}(k)\}.
\ee
With the definition
$\ln n=k{\ln(1/p)\over c}$
we can write:
$$
\Theta(q)=(2k-q)k {\ln(1/p)\over c} +2k-q-q\ln q-2(k-q)\ln(k-q)+{ o}(k) 
$$
We can estimate $\bar Z$ for large $k$ with the saddle point method looking for the maximum
of the function $\Theta(q)+\Phi (q)$. Indeed $\ln \bar Z=\max_{q\in[0,k]}(\Theta(q)+\Phi (q))
+{\cal O}(\ln k)$
We have $\Theta(q)+\Phi(q)= a(q)-b(q)$, where $a(q)$ is a polynomial with degree less or equal $2$, i.e., with
$$
a(q)= (2k-q)k {\ln(1/p)\over c} +2k-q-\b h(k-q)-f(\b)(k^2-q^2)-f(2\b){q(q-1)\over 2}, 
$$
and $b(q)$ the remaining part:
$$
 b(q)=q\ln q+2(k-q)\ln(k-q)+{ o}(k)
$$
By noting that $f(\b)$ is
a concave function so that
$f(\b)-{f(2\b)\over 2}>0$ we have that $a(q)$ is a convex parabola, and so with maximum
 in $0$ or $k$, while $b(q)$ is non negative and $|b(q)|\le 4k\ln k$.
We have that, for sufficiently large $k$,
 the maximum of $a(q)$ is obtained in $q_{max}=k$ if $\tilde h> \tilde h_c$, (see equation (\ref{tildehc}) for the definition of $\tilde h_c$)
and in $q_{max}=0$ if $\tilde h<\tilde h_c$.
By simple calculations we have that these points $k$ and $0$ correspond to the maximum
also for the function $\Theta(q)+\Phi (q)= a(q)-b(q)$ when  $\tilde h> \tilde h_c$
and  $\tilde h<\tilde h_c$, respectively. Indeed in the two different cases it is immediate to
verify that in a neighborhood  of $q_{max}$, i.e., 
  in the intervals
$[k-k^\alpha,k]$ and  $[0,k^\alpha]$ with $\alpha\in(0,1)$, respectively, the following  estimates
hold for the variations of the functions $a$ and $b$: for sufficiently large $k$, there exists
a positive $\e$ such that
$$
|\Delta a(q)|:= |a(q+1)-a(q)|>\e k , \quad |\Delta b(q)|=o(k)
$$
and for $q$ outside these intervals $a(q)<a(q_{max})-\e k^{1+\alpha}$, so that
$a(q)-b(q)\le a(q)<a(q_{max})-\e k^{1+\alpha}<a(q_{max})-b(q_{max})$.

Summarizing we have, for $\tilde h> \tilde h_c$:

\be{barF}
\ln \bar Z= 
k^2 {\ln(1/p)\over c} +k-f(2\b){k(k-1)\over 2}-k\ln k+{ o}(k)
\ee
and for $\tilde h<\tilde h_c$
\be{barF'}
\ln \bar Z=
2k^2 {\ln(1/p)\over c} +2k-\b \tilde hk^2-f(\b)k^2
-2k\ln k+{ o}(k)
\ee
For large $k$ 
we obtain that ${1\over k^2} \ln \bar Z$ is a continuous function with a discontinuous derivative
in $\b$
when $\tilde h=\tilde h_c$, corresponding to a line of a
first order phase transition as discussed in Section \ref{S3}.


\section{The asymptotic self-averaging of $Z$}
The proof of self averaging of $Z$ is a crude calculation based on elementary
arguments. We first evaluate the second moment of $Z$ proving that
asymptotically it behaves like $\bar Z^2$. An upper bound for ${var Z\over \bar Z^2}$
is obtained with a more detailed computation based on the same tools.

\subsection{The second moment of $Z$}
We evaluate the second moment of $Z$:
\be{Z2}
{\mathbb E}(Z^2)={\mathbb E}\Big( \sum_{\s,\t,\s',\t'\in\cX^{(n)}_k} e^{-\b[H(\s,\t)+H(\s',\t')]} \Big).
\ee
By defining, as before,  $q={|\s\cap\t|}$ (and similarly $q'$) we have
\be{hsths't'}
H(\s,\t)+H(\s',\t')=\sum_{\{i,j\}} J_{ij}(\s_i\t_j+\s_j\t_i+\s'_i\t'_j+\s'_j\t'_i)- h (2k-q-q')
\ee
The quantity $\s_i\t_j+\s_j\t_i+\s'_i\t'_j+\s'_j\t'_i$ takes values $0,1,2,3,4$ as in 
the Table 2, where we use the previous notation, i.e.,   $I:=\s \cap \t$, $S:=\s\backslash I$, $T:=\t\backslash I$,
$C:=(\s\cup \t)^c$
and similarly for the sets $I', S', T'$ and $C'$. We also use the notation $SS'$ for the set
$S\cap S'$ and
so on. The table is symmetric due to the symmetry in the exchange  $i\leftrightarrow j$ so
we write only the upper triangle.
\begin{table}
{\scriptsize
\begin{center}
\begin{tabular}{|c|c|c|c|c|c|c|c|c|c|c|c|c|c|c|c|c|}
\hline
$j\backslash i$&SS'&SI'&ST'&SC'&IS'&II'&IT'&IC'&TS'&TI'&TT'&TC'&CS'&CI'&CT'&CC'\\
\hline
SS'&0&1&1&0&1&2&2&1&1&2&2&1&0&1&1&0\\
\hline
SI'&&2&1&0&2&3&2&1&2&3&2&1&1&2&1&0\\
\hline
ST'&&&0&0&2&2&1&1&2&2&1&1&1&1&0&0\\
\hline
SC'&&&&0&1&1&1&1&1&1&1&1&0&0&0&0\\
\hline
IS'&&&&&2&3&3&2&1&2&2&1&0&1&1&0\\
\hline
II'&&&&&&4&3&2&2&3&2&1&1&2&1&0\\
\hline
IT'&&&&&&&2&2&2&2&1&1&1&1&0&0\\
\hline
IC'&&&&&&&&2&1&1&1&1&0&0&0&0\\
\hline
TS'&&&&&&&&&0&1&1&0&0&1&1&0\\
\hline
TI'&&&&&&&&&&2&1&0&1&2&1&0\\
\hline
TT'&&&&&&&&&&&0&0&1&1&0&0\\
\hline
TC'&&&&&&&&&&&&0&0&0&0&0\\
\hline
CS'&&&&&&&&&&&&&0&1&1&0\\
\hline
CI'&&&&&&&&&&&&&&2&1&0\\
\hline
CT'&&&&&&&&&&&&&&&0&0\\
\hline
CC'&&&&&&&&&&&&&&&&0\\
\hline
\end{tabular}
\end{center}
}
\caption{The value of $\s_i\t_j+\s_j\t_i+\s'_i\t'_j+\s'_j\t'_i$ for different $i,j$}
\end{table}
For every $l\in\{1,2,3,4\}$  again  we denote  by ${\cal E}_l$ the set of
unordered pairs $\{i,j\}$ where $\s_i\t_j+\s_j\t_i+\s'_i\t'_j+\s'_j\t'_i=l$.
By the table we have: ${\cal E}_4=II' \times II'$, ${\cal E}_3=(SI' \times \t I')\cup
(IS'\times I\t')\cup (II'\times IT')\cup (II'\times TI')$, and so on.
With these notations we can write for the second moment of $Z$:
$$
{\mathbb E}(Z^2)=\sum_{q=0}^k\sum_{q'=0}^k e^{-\b  h (2k-q-q')}\sum_{\stackrel{\scriptstyle \s,\t:\,|\s\cap\t|=q,}
{\scriptstyle \s',\t':\,|\s'\cap\t'|=q'}} \prod_{\{i,j\}\in{\cal E}_1}{\mathbb E}e^{-\b J_{ij}}\times
$$
\be{Z2_1}
\times\prod_{\{i,j\}\in{\cal E}_2}{\mathbb E}
e^{-2\b J_{ij}}\prod_{\{i,j\}\in{\cal E}_3}{\mathbb E}
e^{-3\b J_{ij}}\prod_{\{i,j\}\in{\cal E}_4}{\mathbb E}
e^{-4\b J_{ij}}
\ee
For shortness we will denote by $g_r $  
the cardinality of the intersection of the different subsets involved
in this table, 
 where the index $r\in\{1,2,...,9\}$ is fixed  in Table 3, e.g  $g_1:={|SS'|}$.
 \begin{table}
\begin{center}
\begin{tabular}{|c|c|c|c|}
\hline
&S'&I'&T'\\
\hline
S&1&2&3\\
\hline
I&4&5&6\\
\hline
T&7&8&9\\
\hline
\end{tabular}
\end{center}
\caption{Index of the intersections}
\end{table}
The cardinalities $g_r$ have the following constraints:
\be{constraints}
g_1+g_2+g_3\le k-q,\quad
g_4+g_5+g_6\le q,\quad
g_7+g_8+g_9\le k-q,
\ee
\be{constraints'}
g_1+g_4+g_7\le k-q',\quad
g_2+g_5+g_8\le q',\quad
g_3+g_6+g_9\le k-q'.
\ee
The cardinality  of the sets given by intersections with $C$
 or $C'$ is obtained by difference:
 \be{CC'}
 g_{S}:=|S C'|=k-q-(g_1+g_2+g_3), \quad   g_{S'}:=|S'C|=k-q'-(g_1+g_4+g_7), 
  \ee
 \be{CC'bis}
 g_{I}:=|I C'|=q-(g_4+g_5+g_6), \quad g_{I'}:=|I' C|=q'-(g_2+g_5+g_8),
 \ee
 \be{CC"}
 g_{T}:=|T C'|=k-q-(g_7+g_8+g_9), \quad  g_{T'}:=|T' C|=k-q'-(g_3+g_6+g_9),
 \ee
 By defining $g=g_1+g_2+...+g_9$ and $M(g_1,...,g_9,n,q,q')$ the multinomial coefficient
 $$
 M(g_1,...,g_9,n,q,q')={n!\over g_1!...g_9! g_{S}!g_{I}!g_{T}!g_{S'}!g_{I'}!g_{T'}!(n-(4k-q-q'-g))!}
$$
 we can write
 \be{Z2_3}
 {\mathbb E}(Z^2)=
 \sum_{q=0}^k\sum_{q'=0}^k \sum_{g=0}^{(2k-q)\wedge(2k-q')}
 \sum_{g_1,...,g_9} M(g_1,...,g_9,n,q,q') e^{\Phi(q)+\Phi(q')+\Psi(q,q',g,g_1,...,g_9)}
 \ee
 where the sum over $g_1,...,g_9$ satisfies the constraints (\ref{constraints}), (\ref{constraints'})
 and  $g_1+g_2+...+g_9=g$ and 
  with $\Phi$ defined in (\ref{phi}) and
  $\Psi$  given by:

 $$
 \Psi={g_5(g_5-1)\over 2}\Big(2f(2\b)-f(4\b)\Big)+{1\over 2}
 \sum_{r=1}^9g_r C_{r}
 $$
 where 
  the coefficients 
 $C_{r}$ are defined as follows:
 \be{C1}
 C_{1}=\Big(2f(\b)-f(2\b)\Big)(g_5+g_6+g_8+g_9)
\ee
\be{C3}
 C_{3}=\Big(2f(\b)-f(2\b)\Big)(g_4+g_5+g_7+g_8)
 \ee
 \be{C7}
  C_{7}=\Big(2f(\b)-f(2\b)\Big)(g_2+g_3+g_5+g_6)
\ee
 \be{C9}
  C_{9}=\Big(2f(\b)-f(2\b)\Big)(g_1+g_2+g_4+g_5)
\ee
 \be{C2}
 C_{2}=\Big(2f(\b)-f(2\b)\Big)(g_4+g_6+g_7+g_9)+
 \Big(f(\b)+f(2\b)-f(3\b)\Big)(g_5+g_8)
 \ee
 \be{C4}
 C_{4}=\Big(2f(\b)-f(2\b)\Big)(g_2+g_3+g_8+g_9)+
 \Big(f(\b)+f(2\b)-f(3\b)\Big)(g_5+g_6)
\ee
 \be{C6}
 C_{6}=\Big(2f(\b)-f(2\b)\Big)(g_1+g_2+g_7+g_8)+
 \Big(f(\b)+f(2\b)-f(3\b)\Big)(g_4+g_5)
 \ee
 \be{C8}
 C_{8}=\Big(2f(\b)-f(2\b)\Big)(g_1+g_3+g_4+g_6)+
 \Big(f(\b)+f(2\b)-f(3\b)\Big)(g_2+g_5)
\ee
\be{C5}
 C_{5}=\Big(2f(\b)-f(2\b)\Big)(g_1+g_3+g_7+g_9)+
 \Big(f(\b)+f(2\b)-f(3\b)\Big)(g_2+g_4+g_6+g_8)
 \ee

We denote  by ${{\cal P}}$ 
the region  of parameters $q,q',g,g_1,...,g_9$ 
defined by the constraints $0\le q\le k$, $0\le q'\le k$, $g=g_1+...+g_9$ and
(\ref{constraints}) and (\ref{constraints'}). 
\bl{L1}
For any $\b\in (0,\infty)$ and for any  $q,q',g_1,...,g_9,g\in {{\cal P}}$ and for $c\le 2$ we have
\be{stima1}
\Psi(q,q',g_1,...,g_9,g)\le \bar\Psi(q,q',g,g_5)
\ee
where
\be{tildephi}
\bar\Psi={g_5(g_5-1)\over 2}\Big(2f(2\b)-f(4\b)\Big)+{1\over 2}\Big(f(\b)+f(2\b)-f(3\b) \Big)((k\wedge g)+g_5)(g-g_5)
\ee
\el
The proof of Lemma \ref{L1} is given in Appendix \ref{A2}.

As far as the entropic term is concerned we can write 
$$
  \sum_{g_1,g_2,g_3,g_4,g_6,g_7,g_8,g_9}{n!\over g_1!...g_9! g_{S}!g_{I}!g_{T}!g_{S'}!g_{I'}!g_{T'}!(n-(4k-q-q'-g))!} \le 
 $$
 \be{entr}
\le \exp\{(4k-q-q'-g)\ln n -\ln\Big( (q-g)!(q'-g)! \Big) -2\ln\Big((k-q-g)!(k-q'-g)!  \Big) +8
 \}=:e^{\bar\Theta_2}
 \ee
 where the sum is under the condition $ g_1+g_2+g_3+g_4+g_6+g_7+g_8+g_9=g-g_5$ and
 with the notation $m!=1$ if $m\le 1$ and 
 where, for the sum of $g_r$ with $r\not= 5$, we used the estimate
 $$
 \sum_{g_1,g_2,g_3,g_4,g_6,g_7,g_8,g_9}{1\over g_1!g_2!g_3!g_4!g_6!g_7!g_8!g_9!}= {8^{g-g_5}\over (g-g_5)!}\le e^8.
 $$
  
 With these estimates we can write
 \be{EZ2_1}
 {\mathbb E}(Z^2)\le
 \sum_{q=0}^k\sum_{q'=0}^k \sum_{g=0}^{(2k-q)\wedge(2k-q')}
 \sum_{g_5=0}^{g\wedge q\wedge q'} e^{\bar\Theta_2(q,q',g,g_5)+\Phi(q)+\Phi(q')+\bar\Psi(q,q',g,g_5)}
 \ee
 \bigskip
 To evaluate this sums  again we look for the maximum of the exponent. Define for notation convenience
 $\Phi_2(q,q')=\Phi(q)+\Phi(q')$.
 \bl{L2}
The maximum of the function  $\bar\Theta_2+\Phi_2+\bar\Psi$ on the parameter region 
defined by the constraints is obtained for $q=q'$ i.e., for  $q=q',g,g_5$ in the three dimensional
polyhedron ${\bar{\cal P}}$ defined by
\be{polyt}
0\le q\le k,\quad g_5\le g\le 2k-q, \quad 0\le g_5\le q
\ee
 and represented in Figure 2.
Moreover $(\bar\Theta_2+\Phi_2+\bar\Psi)(q,q,g,g_5)$ 
reaches its maximum on ${\bf \bar{\cal P}}$ 
in $(k,k,0,0)$ if $\tilde h\ge \tilde h_c$ and in $(0,0,0,0)$ if $\tilde h< \tilde h_c$.
 \el
 
 \begin{figure}
\begin{center}
\includegraphics[width=2in]{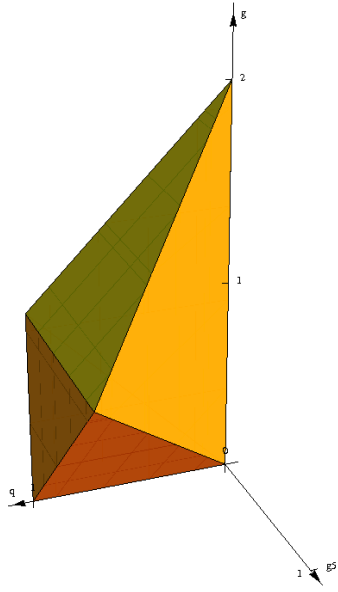}

\caption[]{The polyhedron  ${\bar{\cal P}}$}
\end{center}      
\end{figure}
 
The proof of Lemma \ref{L2} is given again in Appendix \ref{A2}. 
Note that when $g=0$ and $q=q'$ we have the expected relations  $\Theta_2=2\Theta+ {\cal O}({k\over n})$,  and $\Psi=0$.

With these lemmas we immediately obtain
\be{Z2_4}
\ln  {\mathbb E}(Z^2)= \max_{q,g,g_5\in{\bf \bar{\cal P}}}(\bar\Theta_2+\Phi_2+\bar\Psi)(q,q,g,g_5)+
{\cal O}(\ln k)=
2 \ln {\mathbb E}Z+ {\cal O}(\ln k)
 \ee

\subsection{Self averaging}
To evaluate the quantity ${var Z\over({\mathbb E}Z)^2}$
we can write
$$
 ({\mathbb E}Z)^2=
 \sum_{q=0}^k\sum_{q'=0}^k \sum_{g=0}^{(2k-q)\wedge(2k-q')}
 \sum_{g_1,...,g_9} {n!\; e^{\Phi(q)+\Phi(q')}
\over g_1!...g_9! g_{S}!g_{I}!g_{T}!g_{S'}!g_{I'}!g_{T'}!(n-(4k-q-q'-g))!}
$$
and note that, as in the case of the clique number,  the terms corresponding to $g=0$ and $g=1$
are identical in ${\mathbb E}(Z^2)$ and in  $({\mathbb E}Z)^2$, indeed  $\Psi=0$ 
not only for $g=0$ but also in the case $g=1$. Therefore
\be{sa1}
{var Z\over({\mathbb E}Z)^2}\le  {1\over \bar Z^2}
 \sum_{q=0}^k\sum_{q'=0}^k \sum_{g=2}^{(2k-q)\wedge(2k-q')}
 \sum_{g_5=0} ^{g\wedge q\wedge q'} 
 e^{\bar\Theta_2+\Phi_2+\bar\Psi}
\ee
\bl{L3}
The maximum of the function  $(\bar\Theta_2+\Phi_2+\bar\Psi)(q,q',g,g_5)$ on the parameter region 
 ${\bar{\cal P}}$ with the additional constraint $g\ge 2$  is equal to
 $$
 \cases{
-f(2\b)k(k-1)+(2k-2)\ln n -2\ln (k-2)!+o(k)
& if $ \tilde h> \tilde h_c $\cr
 -2\b hk-2f(\b)k^2+(4k-2)\ln n -4\ln (k-2)!+o(k)
& if $  \tilde h< \tilde h_c$\cr
}
 $$
 \el
The proof of this Lemma is analogous to that of Lemma \ref{L2} given in the Appendix \ref{A2}. 

With this Lemma we  conclude the self averaging result. Consider first the case
$\tilde h\ge \tilde h_c$:
$$
{var Z\over({\mathbb E}Z)^2}\le \exp\Big\{-2\Big[ k \ln n +k-f(2\b){k(k-1)\over 2}-k\ln k+{ o}(k)  \Big] +
$$
$$
  -f(2\b)k(k-1)+(2k-2)\ln n -2\ln (k-2)!+o(k) \Big\}\le
$$
$$
\le  \exp\Big\{-2k+2k\ln k -2 \ln n-2(k-2)\ln({k-2\over e})+o(k)\Big\}=
$$
$$
=  \exp\Big\{-2 \ln n+o(k)\Big\}=e^{-2 \ln n+o(k)}
$$
 and using the asymptotic $\ln{k \over k-2}\sim {2\over k}$ we obtain the self averaging in this case.
In the case $\tilde h<\tilde h_c$ the calculation is similar:
$$
{var Z\over({\mathbb E}Z)^2}\le 
\exp\Big\{-2\Big[  2k \ln n +2k-\b \tilde hk^2-f(\b)k^2
-2k\ln k+{ o}(k)
\Big] 
$$
$$-2\b hk-2f(\b)k^2+(4k-2)\ln n -4\ln (k-2)!+o(k)  \Big\}
$$
$$
\le \exp\Big\{-4k+4k\ln k- 2\ln n-4(k-2)\ln({k-2\over e})+{ o}(k)\Big\}\le e^{-2\ln n+o(k)}
$$


\section{Phase transition across $\tilde h_c(T)$}
\label{S3}

By the previous results on the self averaging of  $Z$ 
with 
\be{Zfin}
\bar Z=
\cases{\exp\Big\{ k^2\Big[{\ln 1/p\over c}- {f(2\b)\over 2} \Big] +o(k^2) \Big\}&if $\tilde h> \tilde h_c$\cr
\exp\Big\{ k^2\Big[2{\ln 1/p\over c} - f(\b)-\b\tilde h\Big] +o(k^2) \Big\}&if $\tilde h< \tilde h_c$\cr
}
\ee
we can conclude that the line $\tilde h_c(T)$ represented  in Figure 1
corresponds to a line of a
first order phase transition. Indeed the function $\ln \bar Z$  turns out to be continuous with a discontinuous derivative
in $\b$
when $\tilde h=\tilde h_c$ and we have $-{\partial\over\partial\b}\ln Z=\m_2(H(\s,\t))$
converges almost surely to
\be{Hmedia}
-{\partial\over\partial\b}\ln \bar Z
\cases{
k^2 f'(2\b)+o(k^2)&if $\tilde h> \tilde h_c$\cr
k^2(f'(\b)+\tilde h)+o(k^2)&if $\tilde h< \tilde h_c$\cr
}
\ee

By the convexity property of the function $\ln Z$ we can conclude with standard arguments
that ${\partial\over\partial\b}\lim_{k\to\infty }\ln Z=\lim_{k\to\infty}{\partial\over\partial\b} \ln Z$
and so 
the same result can be obtained by evaluating directly the mean
$\mathbb E\Big[  \m_2(H(\s,\t)) \Big]$

If we look at the model on the
state space of couple of configurations,
with Gibbs measure $\m_2(\s,\t)={1\over Z}e^{-\b H(\s,\t)}$, the two phases correspond to two different mean 
energies.

As far as the second derivative is concerned we have
\be{varH}
{\partial^2\over\partial\b^2}\ln Z=var_{\m_2}(H)=
\cases{
-k^2 f"(2\b)+o(k^2)&if $\tilde h> \tilde h_c$\cr
-k^2f"(\b)+o(k^2)&if $\tilde h< \tilde h_c$\cr
}
\ee
and again the same result can be obtained by evaluating directly the mean on
the $J_{ij}$ of the
variance w.r.t. the pair measure $\m_2$.

In a similar way we can study the overlap $q(\s,\t)$ by computing ${\partial\over\partial \tilde h}\ln Z$.
Indeed $ \m_2( q(\s,\t))=k+{1\over\b k}{\partial\over\partial \tilde h}\ln Z$; we obtain

\be{qmedia}
-{1\over \b k}{\partial\over\partial\tilde h}\ln \bar Z
\cases{
0&if $\tilde h> \tilde h_c$\cr
k&if $\tilde h< \tilde h_c$\cr
}
\ee
so that the two phases have not only different mean energies but also  different mean overlap.
\bigskip


\section{A low temperature phase transition}
\label{S5}

We prove in this section the last claim
of our main theorem.
The proof is divided in three steps.
First, we made a few prelimiray remarks
on the computation of the annealed partition
function $\bar Z$ and we deduce an almost sure
concentration property of the Gibbs measure $\mu_2$.
Second, we translate this concentration
property in a concentration property
of the marginal law $\mu$.
Last, we evaluate the free entropy
of the measure $\pi_\sigma$
for the typical configurations $\sigma$
by proving a last large deviation estimate.

\subsection{Concentration of the Gibbs measure $\m_2$}

An alternative way to compute the annealed partition function
consists in counting the mean number
${\cal N}(q,l_1,l_2)$
of pairs of configurations $(\sigma, \tau)$
with a given overlap $q = |I|:=|\sigma\cap\tau|$,
a given number $l_1$ of missing links inside $I$,
and a given number $l_2$ of missing links between 
$I$ and $T:=\tau\setminus\sigma$, between
$S:=\sigma\setminus\tau$ and $T$,
as well as  
$S$ and $I$.
We get
\be{Zalter}
\bar Z=\sum_{q=0}^k\sum_{l_1=0}^{q(q-1)\over 2}\sum_{l_2=0}^{k^2-q^2}e^{-\b[2l_1+l_2+h(k-q)]}
{\cal N}(q,l_1,l_2).
\ee
To evaluate ${\cal N}(q,l_1,l_2)$ we use the following argument.
Consider the obvious extension of the definition (\ref{H0})
of $H_0(\s,\t)$ to a generic pair $A,B$ of subsets of $V$:
$$H_0(A,B)=\sum_{i,j\in V,\;i\not= j} J_{i,j}
{\bf 1}_A(i){\bf 1}_B(j).$$
For $a\in \{0,...,k\}$ and $l\in \{0,...,ka\}$ let ${\cal A}:=\{A\subset \s^c:\; |A|=a, \hbox { and } H_0(A,\s)=l\}$,
then
\be{EcalA}
\mathbb E|{\cal A}|={n-k\choose a}{ak\choose l} (1-p)^lp^{ak-l}=e^{ak[{\ln 1/p\over c}-I_p({l\over ak})]+o(k^2)}
\ee
where we denote by $I_p$ the large deviation functional
\be{defI}
I_p: x\in[0,1]\mapsto x\ln{x\over 1-p}+(1-x)\ln{1-x\over p}.
\ee
In Appendix \ref{A0} the main properties of this function are recalled; we just mention here
that the function $I_p(x)$ is related to the function $f(\b)$ used in section \ref{S2} by a Legendre transform, indeed we are doing the same computation of $\bar Z$ by using different variables.
A similar computation can be given for subsets of $\s$ and so we can easily conclude that
$${\cal N}(q,l_1,l_2)= n_I(q,l_1)n_{S,T}(q,l_2)$$
with
\be{nI}
n_I(q,l_1)=
{n\choose q}\exp\{ -{q(q-1)\over 2}I_p({ 2l_1\over q(q-1)}) \}
\ee
\be{nST}
\qquad n_{S,T}(q,l_2)={n-q\choose 2(k-q)}
{2(k-q)\choose k-q}
\exp\{ -(k^2-q^2) I_p({l_2\over k^2-q^2}) \}
\ee
and so
\be{calN}
{\cal N}(q,l_1,l_2) = e^{\Theta(q)-{q(q-1)\over 2}I_p(x_1)-(k^2-q^2)I_p(x_2)+o(k)}.
\ee
with $\Theta(q)$ defined in (\ref{theta}), 
$x_1={ 2l_1\over q(q-1)}$ and $x_2={l_2\over k^2-q^2}$,
The sums over $l_1$ and $l_2$ can be written as sums  over $x_1$ and $x_2$ and 
can be estimated with the saddle point method. 
We then obtain for $\bar Z$  the expression given in (\ref{barF}) and (\ref{barF'})
by using the fact that $f$ is the Legendre transform of $-I$.
In the previous sections we used a different approach to estimate $Z$ and $\bar Z$
because this decomposition becomes not easily dealt with
as soon as second moment estimates are involved.
However, the decomposition proposed here is  useful to prove a concentration  
property for the Gibbs measure $\mu_2$.

Three simple remarks will be used in what follows.

\br{r5.2}
As far as equation (\ref{EcalA}) is concerned, we note that the set of values of density of missing links having positive entropy is given by
\be{Xc}
X_c:=\{x\in[0,1]: {\ln 1/p\over c}> I_p({x})\}.
\ee
\er
\br{r5.3}
Note that by the definition of cavity field (\ref{hi}) we have
$$\sum_{i\in T}h_i(S)=\sum_{i\in S}h_i(T)=H_0(S,T)$$ so that
$H_0(S,T)={\sum_{i\in T}h_i(S)+\sum_{i\in S}h_i(T)\over 2}$
\er
\br{r5.4}
The
the density of missing links is typically constant in subsets of a given set.
More precisely
let $A$ and $B$ be a pair of subsets of $V$ with $ |A|=ak, \; |B|=bk$ with $ a,b>0$
and let $\r\in(0,1)$. Then for every $A_1\subset A$ with $|A_1|=a_1k$, $a_1\in (0,a)$
and for every $\r_1\not=\r$ there exists $\d(\r_1)>0$ such that
$$
\mathbb P(H_0(A_1,B)=\r_1a_1bk^2|H_0(A,B)=\r abk^2)<e^{-\d(\r_1) ab k^2}
$$
The proof is an immediate consequence of the convexity of $I_p$ if we note that
$$
\mathbb P(H_0(A_1,B)=\r_1a_1bk^2|H_0(A,B)=\r abk^2)=
$$
$$= \mathbb P(H_0(A_1,B)=\r_1a_1bk^2,\;
H_0(A_2,B)=\r_2 (a-a_1)bk^2) e^{abk^2 I_p(\r)+o(k^2)}\le
$$
$$\le
e^{ -a_1bk^2I_p(\r_1)-(a-a_1)bk^2I_p(\r_2)+abk^2I_p(\r)+o(k^2)}
$$
with $A_2:=A\backslash A_1$ and $\r_2=\r {a\over a-a_1}-\r_1 {a_1\over a-a_1}$.
\bigskip
\er

For $(\sigma, \tau)$ in ${\cal X}_k^{(n)}\times {\cal X}_k^{(n)}$, let us denote
again the overlap by $q(\s,\t)$ and by $H_0(\s,\t)$ the number of missing links between
$\s$ and $\t$.
We define also
\be{barQ}
\bar Q(\s,\t):=\{i\in T:\; h_i(S)\not\in X_c\}\cup \{i\in S:\; h_i(T)\not\in X_c\}
\ee 
that is the set of points in $S\cup T$ with non-typical number of missing links to the other set,
and let $\bar q(\s,\t):=|\bar Q(\s,\t)|$. Even if the energy $H(\s,\t)$ does not depend on this parameter $\bar q(\s,\t)$, as we will show in Section \ref{5.3},  the value of $\bar q(\s,\t)$ is crucial to perform entropy estimates at low temperature when $\tilde h<\tilde h_c$.
We then set, for any $\delta >0$,
\be{Qtip}
{\cal Q}_\d:=
\cases{
   [1-\d,1] & if $\tilde h > \tilde h_{\bar c}$,\cr
   [0,\d] & if $\tilde h < \tilde h_{\bar c}$.
   }
\ee
\be{Htip}
{\cal H}_\d:=
\cases{
    [f'(2\beta)-\d,f'(2\beta)+\d]  & if $\tilde h > \tilde h_{\bar c}$,\cr
    [f'( \beta)-\d, f'(\b)+\d]  & if $\tilde h < \tilde h_{\bar c}$.
}
\ee
and 
\be{barQtip}
\bar {\cal Q}_\d:=
\cases{
   [0,1] & if \quad $\tilde h > \tilde h_{\bar c}$\quad  or \quad $\tilde h < \tilde h_{\bar c}$ and $T>T_c$,\cr
   [1-2\d,1] & if\quad  $\tilde h < \tilde h_{\bar c}$ and $T<T_c$.
   }
\ee

With these intervals of parameters we can define a set of typical pairs of configuration:
\begin{equation}
\Sigma_{2,\delta} := \left\{
  (\sigma,\tau)\in {\cal X}_k^{(n)}\times  {\cal X}_k^{(n)} :\: {q(\sigma,\tau)\over k}\in{\cal Q}_\d,\quad
  {\bar q(\sigma,\tau)\over 2(k-q(\s,\t))}\in\bar{\cal Q}_\d,\quad
  {1\over k^2}H_0(\sigma,\tau)\in{\cal H}_\d\right\}
\end{equation}
where we define $  {\bar q(\sigma,\tau)\over 2(k-q(\s,\t))}=0$ when $q(\s,\t)=k$.
By the self-averaging property of $Z$ we  have
\begin{equation}
{\mathbb E}\left[\mu_2(\Sigma_{2,\delta}^c)\right]
= {\mathbb E}\left[
  {1\over Z}\sum_{(\sigma,\tau)\in\Sigma_{2,\delta}^c}e^{-\beta H(\sigma,\tau)}
\right]\\
= {1\over \bar Z e^{o(k)}}{\mathbb E}\left[
  \sum_{(\sigma,\tau)\in\Sigma_{2,\delta}^c}e^{-\beta H(\sigma,\tau)}
\right].
\end{equation}
Now, using the previous decomposition (\ref{Zalter})  and the fact that 
$I_p$ is strictly convex and more precisely that $I_p''\geq 2$,
we get, with the saddle point method, and using the remarks  \ref{r5.2}, \ref{r5.3}, \ref{r5.4}
that
\begin{equation}
{\mathbb E}\left[\mu_2(\Sigma_{2,\delta}^c)\right]
\leq e^{-C\d^2 k^2}
\end{equation}
for $k$ large enough and a suitable constant $C$.
Indeed
$$
{\mathbb E}\left[\mu_2(\Sigma_{2,\delta}^c)\right]\le
{1\over \bar Z e^{o(k)}}\sum_{q=0,...,k}{}'\sum _{x_1\in[0,1]}\sum _{x_2\in[0,1]}e^{-\b[2l_1+l_2+h(k-q)]}
{\cal N}(q,l_1,l_2)+
$$
\be{mustyp}
+
{1\over \bar Z e^{o(k)}}\sum_{q=0,...,k}\sum _{x_1\in[0,1]}{}'\sum _{x_2\in[0,1]}{}'e^{-\b[2l_1+l_2+h(k-q)]}
{\cal N}(q,l_1,l_2)+e^{-2C\d^2 k^2}
\ee
where the sums  $\sum '$ are with the restriction given by $\s\in (\Sigma_{2,\d})^c$, i.e.,
 $\sum_{q=0,...,k}{}'$ is with the condition $q\not\in{\cal Q}_\d$ and  $\sum_{x_i}'$ are with the condition  $x_1$ is such that $I_p(x_1)+2\b x_1> f(2\b)+C(\d)$
in the case $\tilde h>\tilde h_c$
and in the case  $\tilde h<\tilde h_c$ with the condition $x_2$ is such that $I_p(x_2)+\b x_2> f(\b)+C(\d)$,
with $C(\d)\ge { \d^2\over 2} \min_{x\in [0,1]} I"_p(x) \ge  { \d^2\over 2}2$. Moreover the last term
$e^{-2C\d^2 k^2}$ estimates the mean of the measure of the pairs $(\s,\t)$ such that
${q(\sigma,\tau)\over k}\in{\cal Q}_\d,\quad
  {1\over k^2}H_0(\sigma,\tau)\in{\cal H}_\d$ but $  {\bar q(\sigma,\tau)\over 2(k-q(\s,\t))}\not\in\bar{\cal Q}_\d$. This can be obtained only in the case $\tilde h < \tilde h_{\bar c}$ and $T<T_c$ and in this regime we have that $f'(\b)\not\in X_c$ and the main contribution to $Z$ is given by pairs of disjoint sets
  $S$ and $T$ with $| H_0(S,T)-f'(\b)|<B\d$ for a suitable constant $B$.
  By remarks  \ref{r5.2} and  \ref{r5.3} we can conclude that $S$ (and $T$) can be decomposed
  into two disjoint parts $S=S_1\cup S_2$, with $|S_i|\ge \d k$ with different density of missing link to $T$. The estimate then follows by remark \ref{r5.4}.

We conclude, with Markov inequality and Borel-Cantelli lemma,
that, almost surely,
\begin{equation}
\mu_2(\Sigma_{2,\delta}^c)
\leq e^{-C \d^2 k^2/2}
\end{equation}
for $k$ large enough. 

\subsection{Concentration of the marginal  measure $\m$}

Starting from $\Sigma_{2,\delta}$ we want to define a set $\Sigma_\d$ of {\it typical configurations} in ${\cal X}^{(n)}_k$ with the property that $\s\in\Sigma_\d$ implies that
 $\p_\s$ is concentrated on the configurations
$\t$ such that $(\s,\t)\in\Sigma_{2,\d}$.
 
 To give a precise definition of this set $\Sigma_\d$ we can proceed as follows.
For all $\sigma$ in ${\cal X}^{(n)}_k$ and $\alpha$ , $\bar\a$ and $\rho$ in $[0,1]$, define
\begin{eqnarray}
s_\sigma(\alpha,\bar\a,\rho) &:=& {1\over k^2}\ln\left|\left\{
  \tau\in{\cal X}^{(n)}_k :\: {q(\sigma,\tau)\over k}=\a, \; {\bar q(\s,\t)\over 2(k-\a)}=\bar\a,\; {H_0(\sigma,\tau)\over k^2}= \rho
\right\}\right|\\
\phi_\sigma(\alpha,\bar\a,\rho) &:=& s_\sigma(\alpha,\bar\a,\rho) - \beta(\rho + \tilde h(1-\alpha))\\
\phi_\sigma^* &:=& \max\left\{
  \phi_\sigma(\alpha,\bar\a, \rho) :\:
  \alpha,\bar\a, \rho \in [0,1]
\right\}
\end{eqnarray}
Since the number of possible values of $\alpha$,\,$\bar\a$
and $\rho$ for which $s_\sigma(\alpha,\bar\a,\rho)$ can be non-negative
and finite is only polynomial in $k$, we note that, for all positive $\delta$ 
and large enough $k$, the quantity
$\ln Z_\sigma$ is such that 
\begin{equation}\label{starest}
k^2\phi^*_\sigma \leq \ln Z_\sigma \leq  k^2\phi^*_\sigma + \delta k^2.
\end{equation}

We then set, for any $\delta > 0 $,
\begin{equation}
\Sigma_\delta := \left\{
  \sigma\in {\cal X}^{(n)}_k :\:
 \hbox{ there exist } \alpha\in {\cal Q}_\d,\;\bar\a\in\bar{\cal Q}_\d,\; \rho\in{\cal H}_\d
    \right. \hbox{ such that }
  \phi_\sigma(\alpha,\bar\a,\rho) \geq \phi^*_\sigma -\delta\Big\}
\end{equation}
Note that for a given positive $\delta$,
for $k$ large enough and for all $\sigma$ in $\Sigma_\delta^c$,
\begin{equation}
\pi_\sigma\left(
\{ \tau\in{\cal X}^{(n)}_k :\: {q(\sigma,\tau)\over k}\in{\cal Q}_\d,
{\bar q(\s,\t)\over 2(k-q(\s,\t)}\in \bar{\cal Q}_\d,\; {H_0(\sigma,\tau)\over k^2}\in{\cal H}_\d
\}
\right)
\leq e^{-\delta k^2/2}
\end{equation}
this means that we have a  concentration property of $\p_\s$ implying that for configurations $\s\not\in\Sigma_\d$
the measure $\p_\s$ is concentrated on values of $(\a,\bar\a,\r)$ not in ${\cal Q}_\d\times\bar{\cal Q}_\d\times{\cal H}_\d$.
Now, due to (\ref{measures}) we have, for a given $\delta >0$ and $k$ large enough,
\begin{equation}
\mu_2(\Sigma_{2,\delta}^c)=
\sum_{\s\in \Sigma_\d}\m(\s)\sum_{\t: (\s,\t)\in \Sigma_{2,\delta}^c}\p_\s(\t)+
\sum_{\s\in \Sigma_\d^c}\m(\s)\sum_{\t: (\s,\t)\in \Sigma_{2,\delta}^c}\p_\s(\t)
 \geq \mu(\Sigma_\delta^c) (1 - e^{-\delta k^2/2}).
\end{equation}
We conclude, using the concentration property of $\mu$, that, 
almost surely,
\begin{equation}
\mu(\Sigma_{\delta}^c)
\leq e^{-C\delta^2k^2/3}
\end{equation}
for $k$ large enough. 


\subsection{Conclusion}
\label{5.3}

To estimate the entropy, up to $o(k^2)$
\be{entropy}
S:= -\sum_{\s\in\cX_k}\m(\s)\ln(\m(\s))
=\ln Z-\sum_{\s\in\cX_k}{Z_\s\over Z}\ln Z_\s=
\ln Z+\b\m(F),
\ee
where, for any $\sigma$ in ${\cal X}_k$,
$F(\s):=-{1\over \b}\ln Z_\s$ is the free energy
associated with $\pi_\sigma$,
it is enough to estimate $F(\sigma)$
for all $\sigma$ in $\Sigma_\delta$.
Indeed, 
 $Z$ is self-averaging 
and we estimated $\ln\bar Z$ up to $o(k)$, moreover we have a polynomial uniform upper bound
on $F$ (polynomial in $k$), and an exponential concentration on $\Sigma_\delta$.
This implies that almost
surely, for any $\delta >0$,
\begin{equation}
\left|\mu(F) - \mu(F|\Sigma_\delta)\right| \leq e^{-C\delta^2k^2/4}.
\end{equation}
We will estimate $F(\s)$, i.e., $\ln Z_\s$,  uniformely on $\Sigma_\delta$ in the following cases:
\begin{itemize}
\item[(A)]
  $\tilde h > \tilde h_{\bar c}$,
\item[(B)]
  $\tilde h < \tilde h_{\bar c}$, and $C(\beta) <0$,
  i.e.,  $T<T_{\bar c}$,
\item[(C)]
  $\tilde h < \tilde h_{\bar c}$, and $C(\beta) >0$,
  i.e., $T>T_{\bar c}$.
\end{itemize}
\medskip\par\noindent

For any positive $\delta$,
by definition of $\Sigma_\delta$  for $k$ large large enough 
we have the following estimate for $\ln Z_\s$, for all $\sigma$ in $\Sigma_\delta$,
\be{lnzs}
 \max_{\alpha\in{\cal Q}_\d,\bar\a\in\bar{\cal Q}_\d,\rho\in {\cal H}_\delta} \phi_\sigma (\alpha,\bar\a,\rho)
\leq \phi_\s^*\leq
{1\over k^2}\ln Z_\sigma
\leq  \phi_\s^*+\d\leq
\max_{\alpha\in{\cal Q}_\d,\bar\a\in\bar{\cal Q}_\d,\rho\in {\cal H}_\delta} \phi_\sigma (\alpha,\bar\a,\rho)
+2\delta.
\ee
We estimate $\max_{\alpha\in{\cal Q}_\d,\bar\a\in\bar{\cal Q}_\d,\rho\in {\cal H}_\delta} \phi_\sigma (\alpha,\bar\a,\rho)$ in the 
three different cases.

{\bf Case (A):}
Since $s_\sigma (\alpha,\rho)\geq 0$ and
\begin{equation}
\max_{\alpha\in{\cal Q}_\d,\bar\a\in\bar{\cal Q}_\d,\rho\in {\cal H}_\delta} 
  s_\sigma (\alpha,\bar\a,\rho)
\leq \delta {\ln (1/p) \over c}
\end{equation}
we have
\begin{equation}
-\beta(1+\tilde h)\delta -\beta f'(2\beta)
\leq
{1\over k^2}\ln Z_\sigma
\leq
-\beta f'(2\beta) + (2 + \beta)\delta  
+\delta {\ln (1/p) \over c} .
\end{equation}
by using that $c$ goes to $\bar c$,
we conclude that, almost surely,
\begin{equation}
\lim_{k\rightarrow +\infty} {S\over k^2}
= {\ln (1/p) \over \bar c} - {f(2\beta)\over 2} + \beta f'(2\beta).
\end{equation}

\medskip\par
In cases (B) and (C) we will need
a concentration result on the  numbers of sites $i$ outside $\sigma$
such that $h_i(\sigma) = j + \tilde h k$, i.e., $ g_{j,1}=|{\cal I}_{j,1}|$ (see (\ref{Ij}).
\begin{lemma}\label{nest}
Let
\be{jc}
J_c:=\{j\in\mathbb N :\, {j\over k}\in X_c\}
\ee
with $X_c$ defined in (\ref{Xc}).
With probability 1,
  for any $\delta >0$,
  if $k$ is large enough then,
  for all $\sigma$ in ${\cal X}^{(n)}_k$, for $j\in J_c$ we have:
  \be{gj1}
    \exp\left\{
      k\left(
        -\delta 
        + {\ln(1/p)\over c} - I_p\left({j\over k}\right)
      \right) \right\}
   \le
    g_{j,1} \le
    \exp\left\{
      k\left(
        \delta+
          {\ln(1/p)\over c} - I_p\left({j\over k}\right)
      \right)\right\};
  \ee
 for $j\not\in J_c$ we have 
 \be{gj11}
 g_{j,1}\le e^{k\d}.
 \ee
\end{lemma}
\medskip\par\noindent
{\bf Proof:}
The random variable $ g_{j,1} $ follows a binomial law with
parameters $n-k \leq n$ and 
\begin{equation}
{k\choose j}(1-p)^{j}p^{k-j} = e^{-kI_p(j/k) + o(k)},
\end{equation} 
so that the usual large deviation estimates give
\begin{eqnarray}
&& {\mathbb P}\left(g_{j,1} \geq \exp\left\{
  k\left(
    \delta+\left[
       {\ln(1/p)\over c} - I_p\left({j\over k}\right)
     \right]_+
  \right)
\right\}\right)
\nonumber\\
&& \quad\leq\quad
\exp\left\{
  -k\left(
    \delta+\left[I_p\left({j\over k}\right) 
    - {\ln(1/p)\over c}\right]_+ +o(1)
  \right)
  e^{k(\delta+[\ln(1/p)/c - I_p(j/k)]_+)}
\right\}
\end{eqnarray}
and, if $\ln(1/p)/c \geq I_p(j/k)$, i.e., $j>j_c$
\begin{eqnarray}
{\mathbb P}\left(g_{j,1} \leq \exp\left\{
  k\left(
      - \delta
      +{\ln(1/p)\over c} - I_p\left({j\over k}\right)
    \right)
\right\}\right)
\leq
\exp\left\{
  -e^{k(\ln(1/p)/c - I_p(j/k) + o(1))}
\right\}.
\end{eqnarray}

Since the number of configurations $\sigma$
is not larger than $e^{k^2\ln(1/p)/c}$,
we obtain our result with the Borel-Cantelli lemma.
\hfill$\square$

By Lemma \ref{nest} we can obtain the following results:
\bi
\item[]
In case (B) $T<T_{\bar c}$ i.e.,$f'(\b)\in [0,1]\backslash X_{\bar c}$, there exists a constant $a_2$
such that, almost surely, for all $k$ large enough,
\be{s2}
\max_{\alpha\in{\cal Q}_\d,\bar\a\in\bar{\cal Q}_\d,\rho\in {\cal H}_\delta} 
s_\s(\a,\bar\a, \r)\le a_2\d
\ee
This immediately follows from (\ref{gj11}).
\item[]
In case (C) $T>T_{\bar c}$ i.e., $f'(\b)\in X_{\bar c}$, there exists a constant $a_3$ 
such that, almost surely, for all $k$ large enough,
\be{s3}
\max_{\alpha\in{\cal Q}_\d,\bar\a\in\bar{\cal Q}_\d,\rho\in {\cal H}_\delta} 
s_\s(\a,\bar\a, \r)<{\ln 1/p\over \bar c}-I_p(f'(\b))+a_3\d
\ee
\ei
The proof of this entropy estimates can be found in Appendix \ref{A3}. 
It is absolutely standard but we give it not only for completeness but also
to show that the point of view of  the Fermi statistics is a useful tool.
The main idea is that
in the asymptotics $k\to\infty$, due to the convexity property of $I_p$,
the entropy is essentially due to the sites $i$ with cavity field $h_i(\s)$
such that ${h_i(\s)\over k}\in(f'(\b) +\tilde h-\d, f'(\b) +\tilde h+\d)$  and the number of such sites
is estimated by Lemma \ref{nest}.

With the entropy estimates (\ref{s2}) and  (\ref{s3}) we can easily complete our proof.

\medskip\par\noindent
{\bf Case (B):} 
By equations (\ref{lnzs}) and (\ref{s2})
if $\sigma$ is in $\Sigma_\delta$, almost surely, for all $k$ large enough,
\begin{equation}
-\beta f'(\beta) -\b\d-\b\tilde h 
\leq {1\over k^2}\ln Z_\sigma
\leq -\beta f'(\beta)+\b\d-\b\tilde h(1-\d) + a_2\delta+2\d. 
\end{equation}
We conclude that, almost surely,
\begin{equation}
\lim_{k\rightarrow +\infty} {S\over k^2}
= 2{\ln (1/p) \over \bar c} - {f(\beta)} + \beta f'(\beta).
\end{equation}

\medskip\par\noindent
{\bf Case (C):} Again by equations (\ref{lnzs}) and (\ref{s3})
if $\sigma$ is in $\Sigma_\delta$,
 almost surely, for all $k$ large enough,
\begin{equation}
 {1\over k^2}\ln Z_\sigma
\leq 
-\beta f'(\beta)+\b\d- \b\tilde h(1-\d)+ {\ln(1/p)\over c} - I_p(f'(\beta)) +a_3\delta +2\d. 
\end{equation}
We conclude that, almost surely,
\begin{equation}
\lim_{k\rightarrow +\infty} {S\over k^2}
\geq 2{\ln (1/p) \over \bar c} - {f(\beta)} 
- {\ln(1/p)\over c} + I_p(f'(\beta))  + \beta f'(\beta)
= {\ln (1/p) \over \bar c}.
\end{equation}
where we used that  $f$ is the Legendre transform
of $-I_p$. The opposite estimate
$$\lim_{k\rightarrow +\infty} {S\over k^2}\leq{\ln (1/p) \over \bar c}$$
is trivial.

\appendix
\section{The functions $f(\b)$ and $I_p(x)$}
\label{A0}
We give here some inequalities
for the function $f(\b):=-\ln\big[p+(1-p)e^{-\b}\big]$  defined in the main theorem.

 \begin{figure}
\begin{center}
\includegraphics[width=2in]{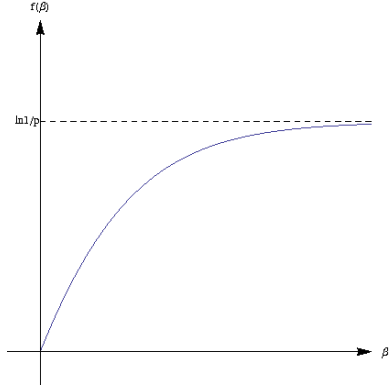}

\caption[]{The function $f(\b)$}
\end{center}      
\end{figure}
This  is a non negative concave function with $f(0)=0$,
$\lim_{\b\to\infty}f(\b)=\ln 1/p=I_p(0)$; from its concavity
we immediately obtain the following estimates:
$$f(\b)>{f(l\b)\over l}\qquad \forall l>1$$
$$B(\b):=f(\b)+f(2\b)-f(3\b)>0$$
since  both the functions 
$F(\b):= f(\b)-{f(l\b)\over l}$ and $B(\b)$ are strictly increasing function
vanishing at zero.
Moreover we have:
$$
f(2\b)+f(3\b)>f(\b)+f(4\b)
$$
which is an immediate consequence of concavity, and 
$$f(\b)+f(3\b)-f(2\b)-{f(4\b)\over 2}>0$$
since $f(\b)+f(3\b)-f(2\b)-{f(4\b)\over 2}=f(\b)-{f(2\b)\over 2}+f(3\b)-{f(2\b)\over 2}-{f(4\b)\over 2}$
again positive by concavity.

For $p\in(0,1)$ and $x\in[0,1]$  consider now the
binomial large deviation functional defined in (\ref{defI}), 
$I_p(x)=x\ln{x\over 1-p}+(1-x)\ln{1-x\over p}$.
This is a convex non negative function with minimum
at $x=1-p$ and $I_p(0)=\ln 1/p, I_p(1)=\ln 1/(1-p)$.
 \begin{figure}
\begin{center}
\includegraphics[width=2in]{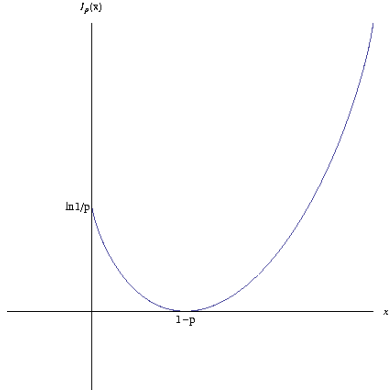}

\caption[]{The function $I_p(x)$ for $ p=2/3$}
\end{center}      
\end{figure}
By recalling the asymptotic behavior for the binomial coefficient:
$$
{L\choose l}\sim (2\pi)^{-1/2}[ x^x(1-x)^{1-x}]^{-L}(x(1-x)L)^{-1/2}
$$
with $x={l\over L}$
 (see for instance \cite{bollobas} pg.4) we immediately obtain
$$
{L\choose l}(1-p)^l p^{L-l}=e^{-LI_p({x})+o(L)}.$$ 

The functions $f(\b)$ and $I_p(x)$ are related by a Legendre transform. Indeed
we have
$$
I_p(x)+\b x\ge f(\b)
$$
where the equality holds only for $x=f'(\b)$.
By evaluating the critical point of the function $I_p(x)+\b x$
we have 
$$
I_p'(x)= \ln {xp\over (1-p)(1-x)}=-\b
$$
ans so the critical point is
$$
x_0={(1-p)e^{-\b}\over p+(1-p)e^{-\b}}= f'(\b)
$$
and this is a minimum due to the convexity of $I_p(x)$.

In particular we have 
$$
I_p"(x)={1\over x(1-x)}\ge 2
$$


\section{Proof of Lemma \ref{L1}}
\label{A1}
 Indeed to prove (\ref{stima1})
we note that, the coefficients $C_{r}$ in
 $$
 \Psi={g_5(g_5-1)\over 2}\Big(2f(2\b)-f(4\b)\Big)+{1\over 2}
 \sum_{r=1}^9g_r C_{r}
 $$
  defined in equations (\ref{C1}), (\ref{C2}), (\ref{C3}), (\ref{C4}), (\ref{C5}), (\ref{C6}), (\ref{C7}), (\ref{C8}), (\ref{C9}), can be estimated by using  the concavity of the function $f(\b)$ so that 
  $0\le 2f(\b)-f(2\b)\le f(\b)+f(2\b)-f(3\b)=:B$.  Indeed  by using the constraints 
  (\ref{constraints}) and (\ref{constraints'}) 
  we can estimate the coefficient:
   $$
 C_{r}\le B(k\wedge g)\quad r\not= 5,\qquad C_5\le B(g-g_5)
 $$
so that
 $$
 \Psi\leq g_5(g_5-1)\Big(f(2\b)-{f(4\b)\over 2}\Big)+{1\over 2}B((k\wedge g)+g_5)(g-g_5)=\bar\Psi
 $$

\section{Proofs of Lemmas \ref{L2} and \ref{L3}}
\label{A2}
{\bf Proof of Lemma \ref{L2}}
We  look now  for the maximum of the function 
$$
\bar\Theta_2+\Phi_2+\bar\Psi= (4k-q-q'-g)\ln n -\ln\Big( (q-g)!(q'-g)! \Big) -2\ln\Big((k-q-g)!(k-q'-g)!  \Big) +C \Big)
 $$
 \be{barest}
+\Phi_2+g_5(g_5-1)\Big(f(2\b)-{f(4\b)\over 2}\Big)+{1\over 2}B((k\wedge g)+g_5)(g-g_5)
\ee
By noting the symmetry of this function in the parameters $q$ and $q'$ and the fact that
the constraints are in the form $g\le (2k-q)\wedge (2k-q')$ and $g_5\le q\wedge q'\wedge g$ , 
we immediately can conclude that the maximum is obtained for $q=q'$. So we have only to
study the function $a(q,g,g_5)-b(q,g,g_5)$ on the polytope ${\bf \bar{\cal P}}$ where
$$
a(q,g,g_5)=-2\b h(k-q)-2f(\b)(k^2-q^2)-f(2\b)q(q-1)+g_5(g_5-1)\Big(f(2\b)-{f(4\b)\over 2}\Big)+
$$
\be{a}
+{1\over 2}B((k\wedge g)+g_5)(g-g_5) +(4k-2q-g)\ln n 
\ee
\be{b}
b(q,g,g_5)=2\ln\Big((q-g)! \Big)+4\ln\Big((k-q-g)! \Big)\Big)+C
\ee
and  ${\cal P}$ (see Figure 2) is defined by the relations:
\be{def_cal_P}
0\le g\le 2k-q,\quad 0\le g_5\le g\wedge q,\quad 0\le q\le k
\ee
We first study the maximum of the function $a$ on $\bar{\cal P}$.
For $g>k$ the hessian of $a$ is given by
\be{hessa1}
\left(
\begin{array}{ccc}
[4f(\b)-2f(2\b)] & 0 &0\\
0 &0&{1\over 2}B\\
0&{1\over 2}B&[2f(\b)-f(4\b)-B]
\end{array}
\right)
\ee
and for $g\le k$ the hessian of $a$ is given by
\be{hessa2}
\left(
\begin{array}{ccc}
[4f(\b)-2f(2\b)] & 0 &0\\
0 &B&{1\over 2}B\\
0&{1\over 2}B&[2f(\b)-f(4\b)-B]
\end{array}
\right)
\ee
Again by the concavity of the function $f(\b)$ in both cases we have a positive eigenvalue
$\l_1=4f(\b)-2f(2\b)$ and two real eigenvalues  with
$\l_2+\l_3>0$ if $g\le k$ and  $\l_2\l_3<0$ if $g>k$ so at least two positive eigenvalues. We can conclude that the maximum
of $a$ is obtained on the edges of $\bar{\cal P}$. 
By studying the function $a({\bf x})$, with ${\bf x}=(q,g,g_5)$, on all the edges we easily check that the maximum actually is 
obtained on the vertices.
To this purpose we used the convexity relations of $f(\b)$ listed in appendix \ref{A0}.
By a direct comparison we obtain that the maximum is obtained 
in the point ${\bf x}_{max}=(k,0,0)$ for $\tilde h> \tilde h_c$ 
and in ${\bf x}_{max}=(0,0,0)$ for $\tilde h< \tilde h_c$
as soon as $f(2\beta) - \frac{1}{2}f(4\beta) < \frac{\ln(1/p)}{c}$.
This inequality holds for all $\beta$ when $c\in(1,2]$, 
while in the case $c>2$ we can simply add
the hypothesis $\b<\bar\b_c$ to conclude.

Fix now $\a\in (0,1)$, 
 in the region
$\bar{\cal P}\cap \{g<k^\a\}$ we have that 
$a({\bf x})-b({\bf x})$ is a decreasing function of $g$ at $q,g_5$ fixed
and  large $k$, and
on the surface $g=g_5$ again is a decreasing function of $g$ for large $k$.
On the other hand we have  for ${\bf x}\in \bar{\cal P}\cap \{g>k^\a\}$ that $a({\bf x})<a({\bf x}_{max})-b({\bf x}_{max})$, so that, as in the discussion of $\bar Z$, by
 noting that  the function $b$ is non-negative, we can conclude that the points ${\bf x}_{max}$
 correspond to maximal values for the function $a({\bf x})-b({\bf x})$.

\section{Proof of equation (\ref{s3})}
\label{A3}

We have to estimate 
\be{nar}
N(\s,\a,\r):=|\{\t\in{\cal X}^{(n)}_k:\, q(\s,\t)=k\a,\, H_0(\s,\t)=k^2\r\}|,
\ee
for $\a\in[0,\d]$ and $\r\in[f'(\b)-\d,f'(\b)+\d]$ with $f'(\b)\in X_c$.
We have  $H_0(\s,\t)=H_0(\s,I)+H_0(\s,T)$ and so  we get
$N(\s,\a,\r)=\sum_{\r'\in[\r-\d,\r+2\d]} N_1(\s,\a,\r') N_2(\s,\a,\r')$
with
$$
N_1(\s,\a,\r')=|\{A\in V\backslash\s:\; |A|=(1-\a)k,\, H_0(\s,A)=k^2(1-\a)\r'\}|,
$$
$$
N_2(\s,\a,\r')=|\{A\in\s:\; |A|=\a k,\, H_0(\s,A)=k^2\a{\r-\r'(1-\a)\over \a}\}|.
$$
The term $N_2$ is easily estimated from above by $2^k=e^{o(k^2)}$.
As far as the term $N_1$ is concerned we can
 use  the notation of the Fermi statistics  and in
particular (\ref{Hen}), to write
\be{N1}
 N_1(\s,\a,\r')= \sum_{\{n_{j,1}\}_{ j=0,1,...,k}:\; \atop
{\sum_{j}n_{j,1}=(1-\a)k, \atop\sum_{j}n_{j,1}j=k^2(1-\a)\r '}} \prod_{j} {g_{j,1}\choose n_{j,1}}
\ee
By using the Stirling formula we can approximate the binomial coefficient
as follows:
\be{stirling}
{g\choose n}= e^{-g {\cal E}({n\over g})+o(k^2)}
\ee
with 
$$
{\cal E}(x):= x\ln x+(1-x)\ln (1-x)
$$
obtaining:
$$
 \sum_{\{n_{j,1}\}_{ j=0,1,...,k}:\; \atop
{\sum_{j}n_{j,1}=(1-\a)k, \atop\sum_{j}n_{j,1}j=k^2(1-\a)\r '}} \prod_{j} {g_{j,1}\choose n_{j,1}}
= \exp\{ \max_{{\bf x}}\sum_{j}[-g_{j,1}({\cal E}(x_j)]+o(k^2)\}
$$
with  ${\bf x}=(x_j)_{j\in \{0,1,...,k\}}$, where $x_j:={n_{j,1}\over g_{j,1}}$, and the maximum is  under the constraints
$\sum_{j}g_{j,1}x_j=(1-\a)k$ and $\sum_{j}g_{j,1}x_jj=k^2(1-\a)\r '$.
With the Lagrange multiplier method and standard computation,  
we can evaluate this maximum by looking at the maximum of the function
\be{Fxlm}
F({\bf x},\l,\m)=\sum_j g_{j,1}\Big[ -{\cal E}(x_j)-(\l+\m j) x_j\Big]
\ee
which is reached in ${\bar{\bf x}}$ with $\bar x_j={1\over 1+e^{\l+\m j}}$ with $\l$ and $\m$ solution of the
equations
\be{lambdamu}
\sum_j g_{j,1}\bar x_j=(1-\a)k \qquad \hbox{ and } \qquad \sum_j g_{j,1}\bar x_j j=k^2(1-\a)\r '.
\ee
In this maximum ${\bar{\bf x}}$ we have
\be{calEbfx}
\sum_{j}[-g_{j,1}({\cal E}(\bar x_j)]=\l (1-\a)k+\m(1-\a)k^2\r '+o(k^2).
\ee
By Lemma  \ref{nest} we have that for $j\in J_c$,  $\bar x_j$ must be exponentially small in $k$
and for any $\d\ge 0$  we have
$$
\sum_{j\in J_c}g_{j,1} x_j= \sum_{j\in J_c}e^{k[{\ln 1/p\over c}-I_p({j\over k})]-\l-k\m {j\over k}+{\cal O}(\d k)}.
$$
Due to the fact that $f'(\b)\in X_c$, this sum is not  exponentially small, i.e., 
$${k[{\ln 1/p\over c}-I_p({j\over k})]-\l-k\m {j\over k}}={\cal O}(\d k)$$
for some $j\in J_c$, and so we can conclude that
$$
\max_{j\in J_c} k[{\ln 1/p\over c}-I_p({j\over k})]-\l-k\m {j\over k}=k[{\ln 1/p\over c}-f(\m)]-\l={\cal O}(\d k)
$$
that is $\l=k[{\ln 1/p\over c}-f(\m)]+{\cal O}(\d k)$ and so,  by (\ref{calEbfx}) that
$$
N_1(\s,\a,\r')\le\exp\Big\{ \{k[{\ln 1/p\over c}-f(\m)]+{\cal O}(\d k)\}(1-\a)k+\m(1-\a)k^2\r ' \Big\}=
$$
$$=\exp\{ k^2(1-\a)[{\ln 1/p\over c}-f(\m)
+\m\r ' ]+{\cal O}(\d k^2)\}
$$
By recalling that  $\r'\in[\r-\d,\r+2\d]=[f'(\b)-2\d,f'(\b)+3\d]$ and 
the Legendre transformation between $f$ and $I_p$ implying that $\m f'(\b)=f(\m)-I_p(f'(\b))$ the proof of (\ref{s2}) and (\ref{s3}) follows straightforward.

\bigskip

{\bf Acknowledgments:}
We thank Antonio Iovanella, coauthor of the numerical parts of this project; the first numerical results he found have been the starting point of our investigations. A special thank is due to Fabio Martinelli for discussions and for the nice and stimulating environment he created in Math Department of Universit\`a di Roma Tre. Thanks to Prasad Tetali for discussions, bibliographic suggestions and encouragement.

\vglue15.truecm

\end{document}